\documentclass[journal,twoside]{IEEEtran}

\usepackage[T1]{fontenc}
\usepackage[latin1]{inputenc}
\usepackage{amssymb,amsfonts,amstext,amsmath}

\usepackage{cite}

   \usepackage[pdftex]{graphicx}
   \DeclareGraphicsExtensions{.eps,.pdf,.jpg}

\usepackage[tight,footnotesize]{subfigure}

\begin{document}
\newcommand{\igen}{j}
\newcommand{\GenSet}{\mathcal{G}}
\newcommand{\GenESet}{\GenSet^{E}}
\newcommand{\GenCSet}{\GenSet^{C}}
\newcommand{\ngen}{J}
\newcommand{\GenBus}{\Gamma_\ibus}
\newcommand{\GenEBus}{\Gamma^E_\ibus}
\newcommand{\GenCBus}{\Gamma^C_\ibus}

\newcommand{\nload}{N_d}
\newcommand{\defcost}{h}

\newcommand{\SpaceState}{\Xi}

\newcommand{\ibus}{i}
\newcommand{\BusSet}{\mathcal{N}}
\newcommand{\nbus}{N}

\newcommand{\ilin}{k}
\newcommand{\LinSet}{\mathcal{K}}
\newcommand{\LinESet}{\LinSet^{E}}
\newcommand{\LinCSet}{\LinSet^{C}}
\newcommand{\nlin}{M}
\newcommand{\LinBus}{\Phi_\ibus}
\newcommand{\LinEBus}{\Phi^E_\ibus}
\newcommand{\LinCBus}{\Phi^C_\ibus}

\title{Reliability-Constrained Power System Expansion Planning: A Stochastic Risk-Averse Optimization Approach
}

\author{Luiz~Carlos~da~Costa~Jr., \and
        Fernanda~Souza~Thom\'{e}, \and
        Joaquim~Dias~Garcia, \and
        Mario~V.~F.~Pereira
}




\maketitle

\begin{abstract}
This work presents a methodology to incorporate reliability constraints in the optimal power systems expansion planning problem. Besides LOLP and EPNS, traditionally used in power systems, this work proposes the use of the risk measures VaR (Value-at-Risk) and CVaR (Conditional Value-at-Risk), widely used in financial markets. The explicit consideration of reliability constraints in the planning problem can be an extremely hard task and, in order to minimize computational effort, this work applies the Benders decomposition technique splitting the expansion planning problem into an investment problem and two sub-problems to evaluate the system's operation cost and the reliability index. The operation sub-problem is solved by Stochastic Dual Dynamic Programming (SDDP) and the reliability sub-problem by Monte Carlo simulation. The proposed methodology is applied to the real problem of optimal expansion planning of the Bolivian power system.
\end{abstract}
\IEEEpeerreviewmaketitle

\section{Introduction}

\IEEEPARstart{T}{he} goal of power systems expansion planning (SEP) is to determine necessary changes in the system due to load growth, new technologies and policy related constraints. In this sense, new generators should be constructed with the goal of satisfying the new system's needs and the planning process decisions are associated to the selection of the best set of equipment (such as generators, transmission lines, transformers, etc.) to accomplish this task. This decision process gives origin to a complex optimization problem, where the objective is to plan the future power system minimizing the investment and operation costs subject to a pre-established set of constraints. Cases of SEP are generation expansion planning (GEP) \cite{zhu1997review,oree2017generation} and transmission expansion planning (TEP) \cite{lumbreras2016new,hemmati2013comprehensive} that can be combined and integrated \cite{pozo2013three,moreira2017reliable} generally called SEP or GTEP.

This planning process constitutes an extremely complex problem that cannot be solved without simplifications. Depending on the goals of the planner different aspects and details of power systems can be considered in general SEP. Many possible constraints are described in \cite{jenkins2017enhanced}, some specific aspects include: carbon capture and storage \cite{saboori2016considering}; unit commitment in the operation \cite{hua2018representing}; flexible demand and electric vehicles \cite{ramirez2016co}; aggressive wind power penetration \cite{zhan2017generation}. A common aspect in most SEP models is the representation of uncertainties, although each model typically focuses on sources of randomness, like renewable energy and load \cite{jabr2013robust}, outages or contingencies \cite{roh2009market,moreira2018five}. Frequently used frameworks to deal with uncertainties in SEP are Stochastic Optimization \cite{jirutitijaroen2008reliability,lopez2007generation} and Robust Optimization \cite{jabr2013robust,mejia2014adjustable}, both of which can also be combined \cite{baringo2018stochastic}. 

Many techniques were proposed to solve the large-scale problems that arise from SEP modelling. Heuristics like Particle Swarm Optimization and GRASP were proposed in \cite{saboori2016considering} and \cite{binato2001greedy}. Many decomposition methods were presented due to natural scenario-wise and/or stage-wise structure: \cite{liu2018multistage} \cite{munoz2015scalable} and \cite{wu2014stochastic} apply progressive hedging, Dantzig Wolfe decomposition was applied in \cite{singh2009dantzig} and Benders decomposition, perhaps the most used one, was applied in \cite{binato2001new,gorenstin1993power,campodonico2003expansion,roh2009market,jirutitijaroen2008reliability,baringo2013risk,jabr2013robust,dehghan2016reliability}.

Reliability is an important aspect of power systems that can be considered in SEP. The simplest way to taking reliability into account is adopting a hierarchical approach, in which the expansion plan is at first elaborated under the economic focus (first step), in other words, aiming the minimization of the investment costs plus the cost of load supply (operation cost). Hereafter, the necessary additional investments to meet a minimum criterion of security (reliability reinforcements) are evaluated (second step).

Among the first approaches for the solution of the SEP with reliability constraints is the work \cite{Cot1979RelOptGenPlan} where linear approximations of a reliability function were obtained from a non-linear formulation. Afterwards, \cite{Blo1983GenCapExpBenders} presented a similar model, but able to generate Benders cuts from a probabilistic simulation model. However, this last model was non-convex which led to convergence problems. Robust optimization was used to induce expansion plans with better reliability indexes in \cite{moreira2017reliable}. In \cite{jirutitijaroen2008reliability}, reliability requirements were modeled by including a cost of load loss in the objective function.

In general, SEP can be formulated as a minimum total cost function (investment and operation) subject to the operation and reliability constraints that depend directly upon the investment decisions. This formulation contains a very opportune structure for the application of decomposition techniques and this characteristic was first explored by \cite{Per1985SensGTExpPlan}. Thereafter, \cite{Oli1987capexp} presented a model for load-peak capacity expansion considering EPNS constraints, but its application was restricted to the second step problem, in other words, the evaluation of additional investments to attend the reliability requisites. Reliability constraints were also applied in \cite{sirikum2006power,park2000improved}, but the problem was solved by genetic algorithm without optimality proofs. More recently, \cite{roh2009market,dehghan2016reliability} considered explicit constraints on reliability indexes via Monte Carlo simulation and Benders decomposition with feasibility cuts. \cite{rashidaee2018linear} considered a Loss Of Load Probability (LOLP) constraint, originally modeled as a Mixed Integer Nonlinear Program (MINLP) and approximated by a Mixed Integer Linear Program (MIP). 

This work proposes an integrated methodology for the solution of the SEP, where the economic and reliability analyses are carried out on an integrated problem. Therefore, it is possible to assess the benefit of each project both in terms of reduction of the operative cost as well as the increase in the overall system reliability. Additionally, we compare the obtained results from the proposed integrated approach with the simplification made by using the hierarchical approach for the solution of the SEP.

We use the Benders decomposition technique to split the original problem in investment, operation and reliability modules. This partition allows each subproblem to be solved by a specialized algorithm, for instance, the investment master problem (a mixed integer programming problem) is solved by \emph{Branch-and-Bound} (B\&B), the operation subproblem by \emph{Stochastic Dual Dynamic Programming} (SDDP) \cite{PerPin1991MultiStageStochasticOptimization} and the reliability subproblem by \emph{Monte Carlo simulation} (MC). 
Instead of focusing on classical measures like LOLP and EPNS \cite{dehghan2016reliability} we propose the use of the two risk measures $\text{VaR}_{\alpha}$ and $\text{CVaR}_{\alpha}$, commonly used in the financial area, but also applied in different contexts in SEP \cite{baringo2013risk,maceira2014application,tohidi2018modified}.

The main contributions of this paper are: a general description of SEP that can be solved by Benders decomposition; application of alternative metrics to power system reliability constraints; description of the reliability constraint as linear programming deterministic problem equivalent to the MC version of the constraint evaluation; combination of the SDDP method as a technique to solve the large-scale convex operation problem.

The remainder of this work is organized as follows. In the next section, \ref{SEP_base}, a generic formulation of reliability constrained system expansion problem is presented. Section \ref{Reliability} describes reliability metrics and their corresponding formulations as optimization problems. In section \ref{Decomposition}, we describe the decomposition procedure that combines the plan choice, cost and reliability evaluation, after that, we present case studied on section \ref{Study}. Finally, conclusions are drawn in section \ref{Conclusion}.

\section{Power Systems Expansion Planning} \label{SEP_base}
The reliability constrained SEP can be formulated as the following mixed integer programming problem:
\begin{subequations}
\label{eq:ProbPlan_General}
\begin{alignat}{2}
	\text{Min} \quad
	& I(x) + O(x) \\
	\text{s.t.:} \quad
	& R(x) \leq \bar{R} \label{eq:ProbPlan_General_Cnf} \\
	& x \in \mathcal{X}
\end{alignat}
\end{subequations}
where $x$ is the vector of investment decisions, $I(x)$ the investment cost, $O(x)$ the operation cost and $R(x)$ the risk measure functions, $\bar{R}$ the reliability criteria and $\mathcal{X}$ the set of decisions that meet  planning constraints.

\subsection{Economic Planning (EP)}
\label{cap:PE}
The first step of an hierarchical planning process is a simplification of problem \eqref{eq:ProbPlan_General}, by disregarding the explicit representation of the reliability aspects and the risk constraints \eqref{eq:ProbPlan_General_Cnf}, as shown in problem \eqref{eq:ProbPlan_InvOpe}.
\begin{subequations}
\label{eq:ProbPlan_InvOpe}
\begin{alignat}{2}
	\text{Min} \quad
	& I(x) + O(x) \\
	\text{s.t.:} \quad
	& x \in \mathcal{X}
\end{alignat}
\end{subequations}

In this sense, the explicit representation of the uncertainty associated to the state of each generating unit is, in general, simplified by applying a reduction factor in every power plants' capacity corresponding to their average availability rate. This approximation avoids the representation of all power system states and the SEP of generation systems. A simplified formulation of EP is given by:
\begin{subequations}
\label{eq:ProbPlan_InvOpeFull}
\begin{alignat}{3}
	\text{Min} \quad
	& \sum_{\igen \in \GenSet} c_\igen x_{\igen} + \sum_{\igen \in \GenSet} d_\igen g_{\igen} + \defcost r & \label{eq:ProbPlan_InvOpeFull_fobj} \\
	\text{s.t.:} \quad
	& \sum_{\igen \in \GenSet} g_{\igen} + r = D & \label{eq:ProbPlan_InvOpeFull_Demanda} \\
	& g_{\igen} \leq \tilde{g}_{\igen} x_{\igen} & \igen \in \GenSet \label{eq:ProbPlan_InvOpeFull_LimGC} \\
 	& x \in \mathcal{X} \label{eq:ProbPlan_InvOpeFull_Inv}
\end{alignat}
\end{subequations}

\noindent where $\GenSet$ is the set generators, $c_\igen$ and $d_\igen$ represent the investment and variable operation costs of generating unit $\igen$, $D$ represents the total system demand, $\defcost$ is the load shedding unit cost, $r$ is a variable that represents the load shedding, $\tilde{g}_{\igen}$ is the available capacity of generator calculated in terms of average availability as $\tilde{g}_\igen = (1 - p_\igen) \times \bar{g}_\igen$, $p_\igen$ is the average failure rate and $\bar{g}_\igen$ is the installed capacity of $\igen$.

The objective function \eqref{eq:ProbPlan_InvOpeFull_fobj} is the minimization of the total cost (investment and operation), subject to demand supply in each time step \eqref{eq:ProbPlan_InvOpeFull_Demanda}, limits on generation \eqref{eq:ProbPlan_InvOpeFull_LimGC}, and investment constraints \eqref{eq:ProbPlan_InvOpeFull_Inv}.

The operation problem is usually much more complex that the illustrative model presented above. As mentioned in the introductory section, many details might be taken into consideration such as: demand response, storage equipment and time coupling, uncertainty, competition, unit commitment, energy network representation and many others. 

\begin{figure}[!t]
	\centering
	\includegraphics[width=1.0\columnwidth]{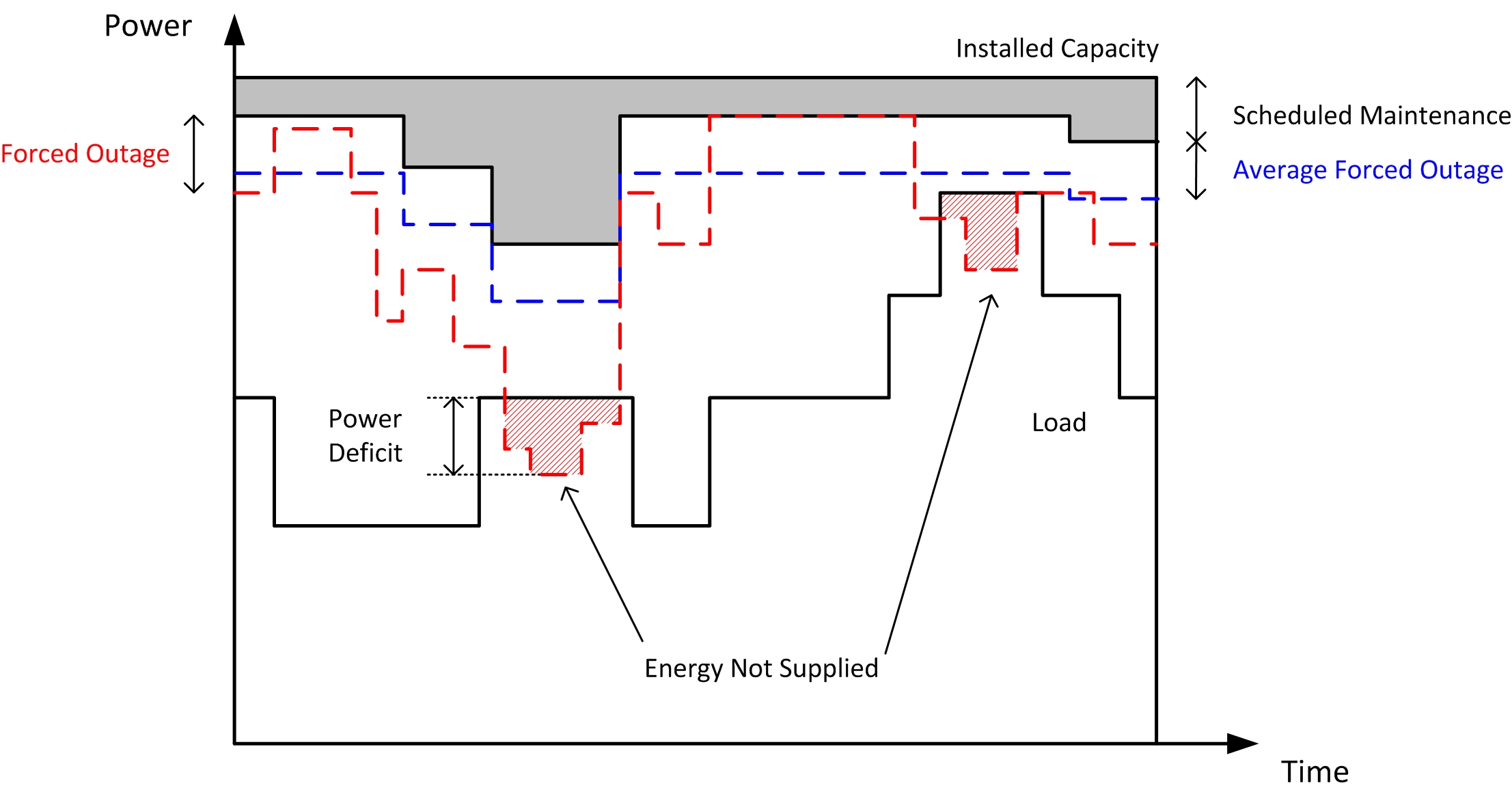}
 	\caption{Random behavior of total system capacity}
 	\label{fig:suprhist}
\end{figure}

Although this model has been frequently used, the use of a generation capacity based on average availability may not be enough to capture the true exposure/risk of load shedding events, as illustrated in Figure \ref{fig:suprhist}. Even if the system can meet the load on average, it might exist one or more states of failure in which the remaining capacity is no sufficient to supply the demand (hatched regions). Hereof it might be necessary to model the system in probabilistic terms to better represent the reliability aspects and this requires modeling the operating state of all generators and, consequently, the total capacity of the system as a random variable (r.v.).

\section{Reliability Analysis of Power Systems} \label{Reliability}
A power system is basically composed of elements such as generators, transmission lines, transformers and load and each element may be in a state among a set of possible states. For example, the operating state of a generating unit can be either (a) $0$ if the equipment is not working or (b) $1$ if it is working. Other elements, such as demand and combined cycle plants, may require a multi-state representation.

A state of a power system with $\ngen$ generators is represented by the random vector $\xi = (\xi_{1}, \xi_{2}, \ldots, \xi_{\ngen})$ where $\xi_{j}$ is a r.v. representing the state of the $\igen$-th generator. Let $S$ represent the set of states of the power system, given by the combination of all possible states of each element and e each state is denoted by $\xi^s, \ s \in S$.
For each state of generator $\igen$ there is an associated probability of occurrence $p_\igen = P(\xi_\igen)$ and, once the state of each generator for the system state $s \in S$ is given, it is possible to calculate the probability of this system state as $p^s = P(\xi^s)$.

The total system capacity is denoted by the r.v. $\bar{G} = \sum_{\igen=1}^{\ngen}\xi_\igen \bar{g}_\igen$ and the load shedding is given by the r.v. $R = \max(D-\bar{G},0)$, which is the insufficient generation capacity. Since $\xi$ has finite support distribution, each realization of system state $\xi^s$ is associated to a total capacity of $\bar{G}^s$ and a respective load shedding $R^s$.

The performance of a given investment plan $x$ can be measured with risk indexes based, in general, on the probability distribution of the load shedding. Since the objective of the SEP is to determine which and when generators should be constructed, the probability distribution of $P(R)$ must also be a function of the investment decision vector $x$. Therefore, the objective of the reliability-constrained SEP becomes to find the plan that minimizes the investment and operation costs and present a ``controlled'' load shedding distribution function in the sense that its associated reliability index meets some given pre-established planning criterion.

The next section presents the reliability indexes LOLP and EPNS, traditionally used in power systems, and then introduces in the context of power system reliability the indexes $\text{VaR}_{\alpha}$ and $\text{CVaR}_{\alpha}$, frequently used in the financial sector. Moreover, the calculation of these indexes will be formulated as optimization problems aiming their incorporation into the SEP.

\subsection{Typical Reliability Indexes}
\begin{figure}[!t]
	\centering
	\subfigure[LOLP]{
		\includegraphics[width=.45\columnwidth]{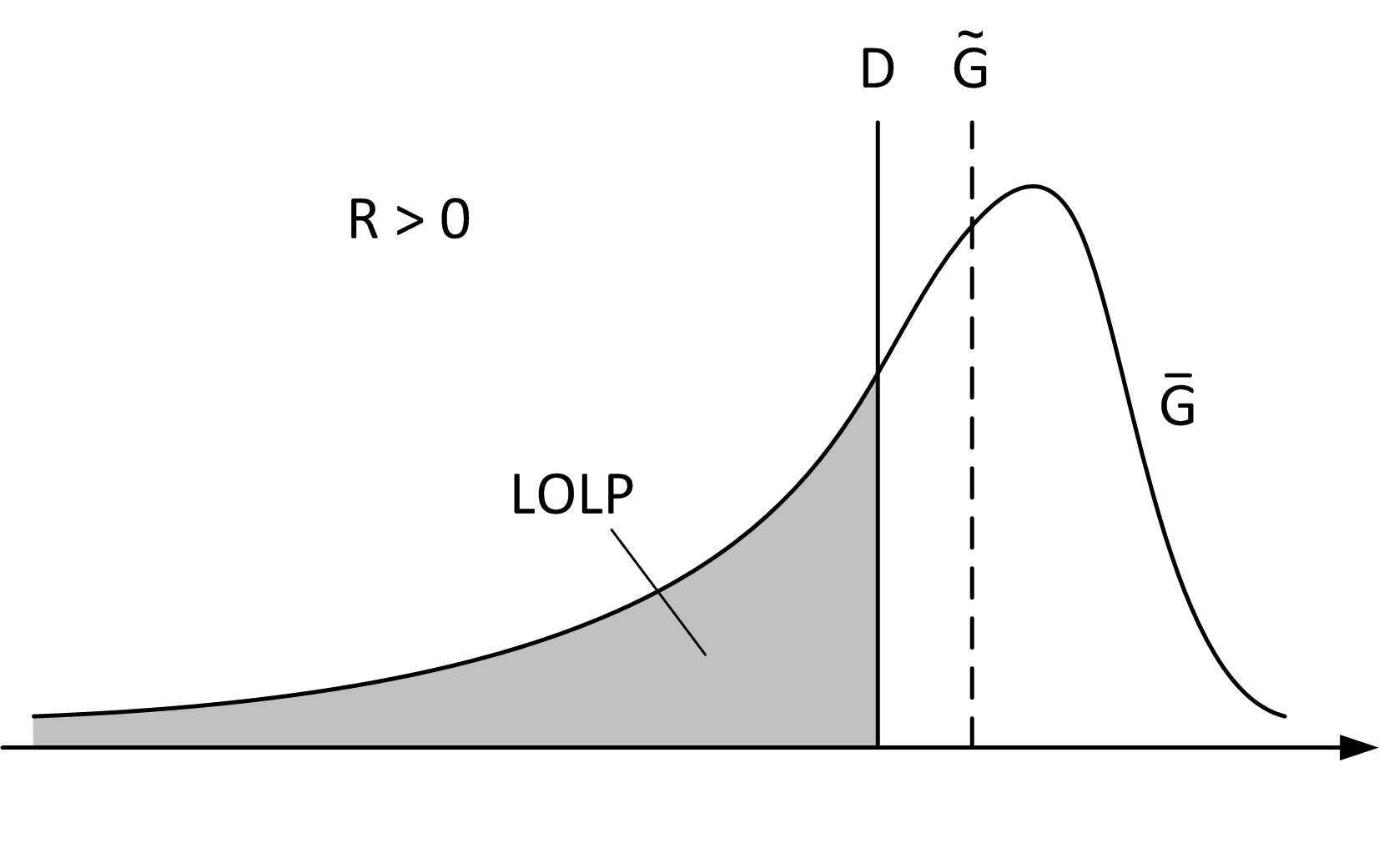}
		\label{fig:dist-lolp}
	}
	\subfigure[EPNS]{
		\includegraphics[width=.45\columnwidth]{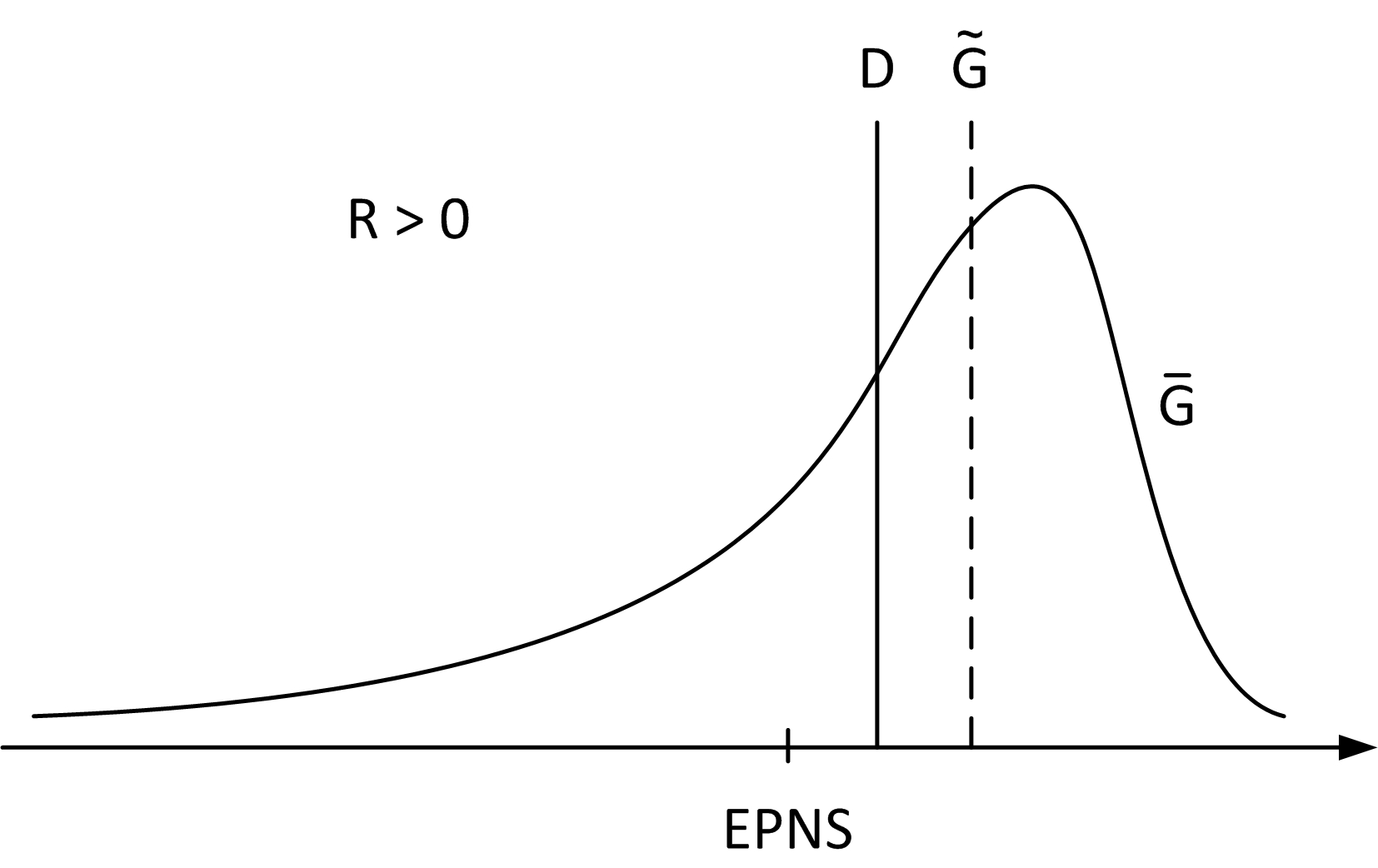}
		\label{fig:dist-epns}
	}
 	\caption{Typical Reliability Indexes}
\end{figure}

\subsubsection{LOLP}
The most straightforward approach to measure the risk of failure of supplying the power system load is to assess the number of insufficient states among all of them. The LOLP (\emph{Loss Of Load Probability}) is the probability of load shedding, as illustrated in Figure \ref{fig:dist-lolp}, and is given by
\begin{equation}
	\text{LOLP} = P(R>0) = \sum_{s \in \Omega} p^s
\end{equation}
where $\Omega = \{ s \in S | r^s > 0 \}$.

The calculation of LOLP can be formulated as a (possibly large) mixed integer linear programming problem, using an integer variable $\phi^s$, for each state $s$ to indicate whether it leads to a load shedding. Then, the LOLP is calculated as the average of these indicator variables weighted by the state probability. By explicitly incorporating this formulation into the problem, the LOLP-constrained EPP can be defined as:
\begin{subequations}
\label{eq:ProbPlan_InvOpeCnf_LOLP}
\begin{alignat}{5}
	\text{Min} \quad
	    & \sum_{\igen \in \GenSet} c_\igen x_\igen + O(x) & \\
	\text{s.t.:} \quad
	    & \sum_{s \in S} p^s \phi^s \leq \overline{\text{LOLP}} \label{eq:ProbPlan_InvOpeCnf_LOLP_LimiteLOLP} \\
	    & r^s \geq D - \sum_{\igen \in \GenSet} \xi_\igen^s \bar{g}_\igen x_\igen     & \qquad \forall s \in S \label{eq:ProbPlan_InvOpeCnf_LOLP_corte} \\
	    & \phi^s \geq \frac{1}{D} r^s
	        & \qquad \forall s \in S \label{eq:ProbPlan_InvOpeCnf_LOLP_ind} \\
	    & r^s \geq 0
	        & \qquad \forall s \in S \\
	    & \phi^s \in \{ 0,1 \}
	        & \qquad \forall s \in S \\
	    & x \in \mathcal{X}
\end{alignat}
\end{subequations}
where $\overline{\text{LOLP}} \in [0,1]$ is the accepted level of reliability adopted as the planning criterion. Constraint \eqref{eq:ProbPlan_InvOpeCnf_LOLP_corte} relates the load shedding to the system's capability of load supply in each state and constraint \eqref{eq:ProbPlan_InvOpeCnf_LOLP_ind} ensures that the indicator variable $\phi_s$ will be equal to $1$ for the states with load shedding. Note that constraint \eqref{eq:ProbPlan_InvOpeCnf_LOLP_LimiteLOLP} limits the value of $\text{LOLP}$ and, consequently, restricts the set of possible investment plans. 

One characteristic of LOLP, which is also the most common source of criticism, is that the depth of the load shedding is disregarded since ``bad'' states are equally labeled with $\phi^s = 1$, independent of the amount of load being shed. This can mislead the SEP to find investment plans with small probability of failure but exposed to states with high level of severity.

\subsubsection{EPNS}
The \emph{Expected Power Not Supplied} is the average value of the load shedding of all system states, as illustrated in Figure \ref{fig:dist-epns}. It can be defined as
\begin{equation}
	\text{EPNS} = E[R] = \sum_{s \in S} p^s r^s
\end{equation}

The EPNS can also be formulated as an optimization problem and explicitly incorporated into the SEP, resulting in the EPNS-constrained SEP model \eqref{eq:ProbPlan_InvOpeCnf_EPNS}.
\begin{subequations}
\label{eq:ProbPlan_InvOpeCnf_EPNS}
\begin{alignat}{5}
	\text{Min} \quad
	& \sum_{\igen \in \GenSet} c_\igen x_\igen + O(x) & \\
	\text{s.t.:} \quad
	& \sum_{s \in S} p^s r^s \leq \overline{\text{EPNS}}\\
	& r^s \geq D - \sum_{\igen \in \GenSet} \xi_\igen^s \bar{g}_\igen x_\igen & \qquad \forall s \in S \\
	& r^s \geq 0 & \qquad \forall s \in S \\
	& x \in \mathcal{X}
\end{alignat}
\end{subequations}
where $\overline{\text{EPNS}}$ is the pre-established planning criterion to EPNS.

Although the EPNS captures the severity of the load shedding in average terms, this index also considers all the ``good'' states (i.e. without load shedding) and, thus, results on a ``diluted'' index, not reflecting the real exposure of the system to states with failure.

\subsection{Risk Measures Used in Financial Area}
It would be interesting if the reliability index could capture both the characteristics of LOLP and EPNS, i.e., both the number of states with load shedding and the severity of these states. To accomplish this task we propose the risk measures used in portfolio optimization problems in financial sector. In the next sections, the risk measures $\text{VaR}_{\alpha}$ and $\text{CVaR}_{\alpha}$ will be presented and introduced in context of reliability analysis for power systems expansion planning.

\subsubsection{$\text{VaR}_{\alpha}$}
\begin{figure}[!t]
	\centering
	\subfigure[$\text{VaR}_{\alpha}$]{
		\includegraphics[width=.45\columnwidth]{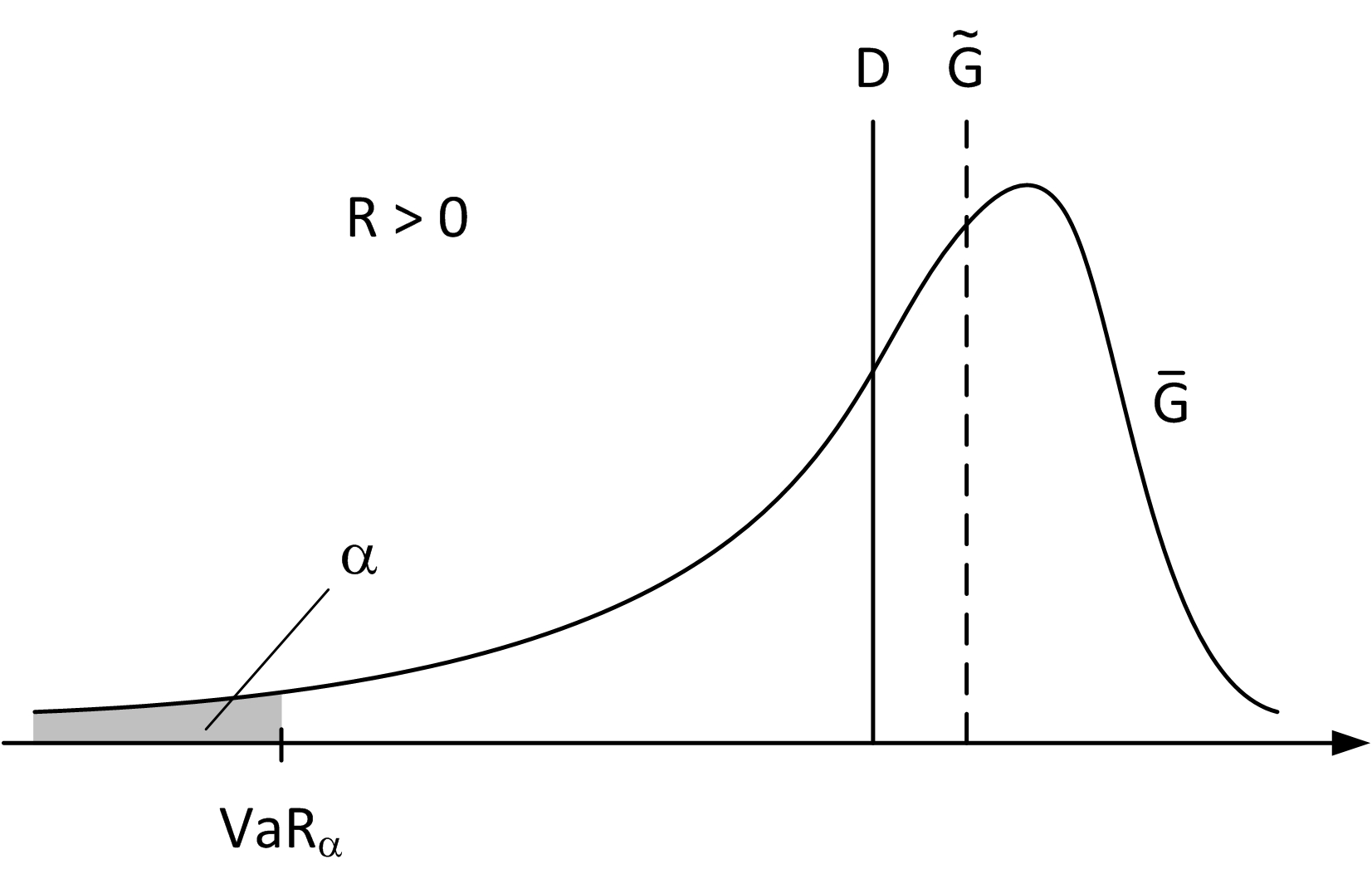}
		\label{fig:dist-var}
	}
	\subfigure[$\text{CVaR}_{\alpha}$]{
		\includegraphics[width=.45\columnwidth]{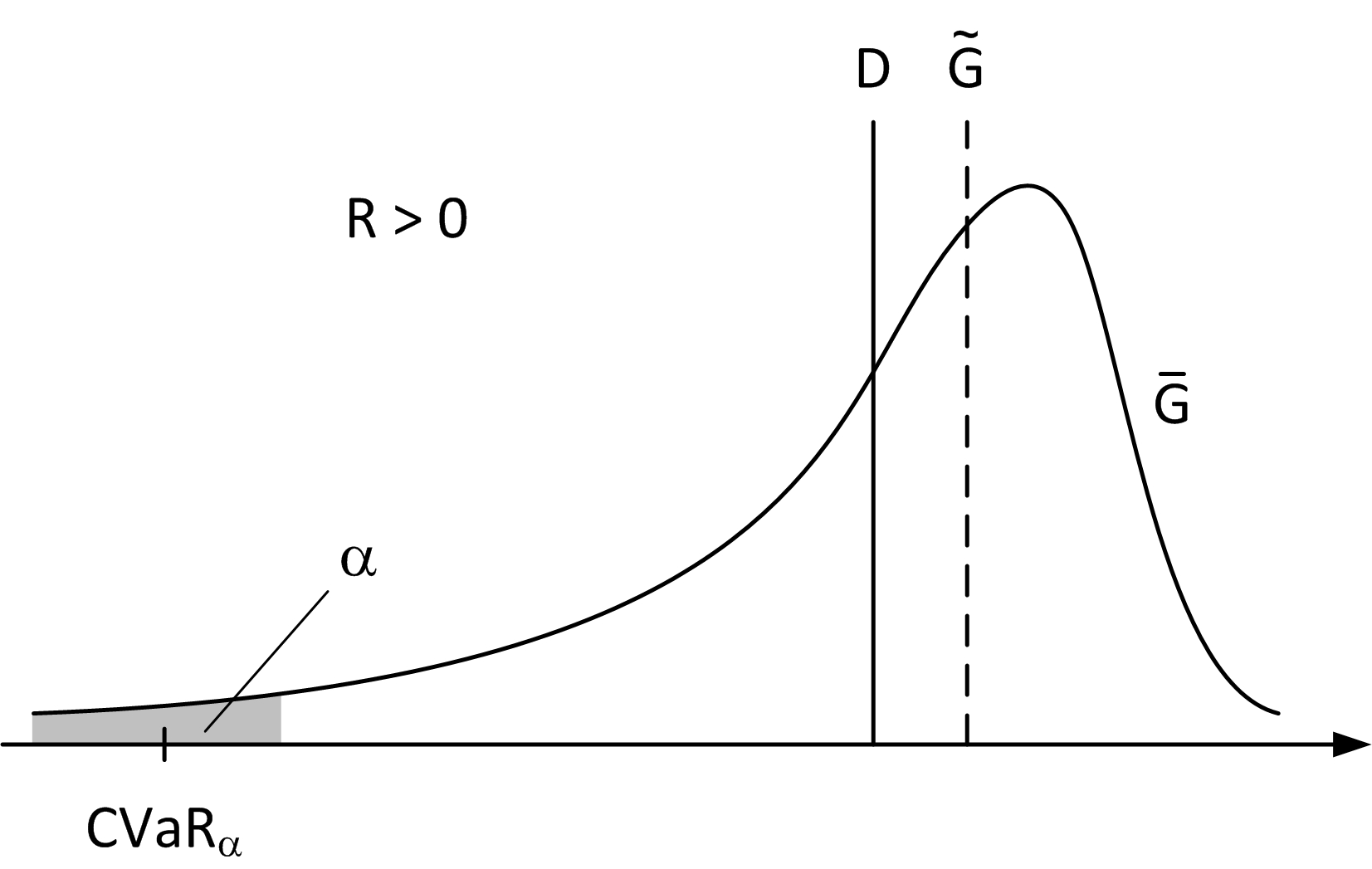}
		\label{fig:dist-cvar}
	}
 	\caption{Financial Risk Measures}
\end{figure}

The \emph{Value-at-Risk} \cite{Jor2000VaR} is a risk index that aims to measure the lowest load shedding associated to a probability of occurrence $\alpha$ or, similarly, the maximum load shedding within a specified level of confidence $1-\alpha$, as illustrated in Figure \ref{fig:dist-var}. For example, $\text{VaR}_{5\%}$ answers the question ``what is the maximum possible load shedding considering the 95\% \emph{best} states''?

As in the case of LOLP and EPNS, it is possible to define the SEP using $R(x) = \text{VaR}_\alpha(x)$ as the reliability index. Its formulation can be explicitly incorporated into the problem  \eqref{eq:ProbPlan_General}, resulting in problem \eqref{eq:ProbPlan_InvOpeCnf_VaR}.
\begin{subequations}
\label{eq:ProbPlan_InvOpeCnf_VaR}
\begin{alignat}{5}
	\text{Min} \quad
	& \sum_{\igen \in \GenSet} c_\igen x_\igen + O(x) & \\
	\text{s.t.:} \quad
	& r^s - D \phi^s \leq \overline{\text{VaR}} & \qquad \forall s \in S \label {eq:ProbPlan_InvOpeCnf_VaR_maxvar} \\
	& \sum_{s \in S} p^s \phi^s \leq \alpha \label{eq:ProbPlan_InvOpeCnf_VaR_sumalpha} \\
	& r^s \geq D - \sum_{\igen \in \GenSet} \xi_\igen^s \bar{g}_\igen x_\igen & \qquad \forall s \in S \\
	& r^s \geq 0 & \qquad \forall s \in S \\
	& \phi^s \in \{ 0,1 \} & \qquad \forall s \in S \\
	& x \in \mathcal{X}
\end{alignat}
\end{subequations}
where $\overline{\text{VaR}}_{\alpha}$ is the limit determined by the system planner. Note that, as in the case of LOLP, an integer variable for each state is also required, what makes its representation more difficult.

Rearranging the constraint \eqref{eq:ProbPlan_InvOpeCnf_VaR_maxvar}, we have $\phi^s \geq D^{-1} \left( r^s - \overline{\text{VaR}} \right)$ which indicates that, when the load shedding $r^s$ exceeds the limit $\overline{\text{VaR}}$ the variable $\phi^s$ should assume value $1$, characterizing the states in the tail of the probability distribution function. Even though it incorporates the parameter $\alpha$ that allows this index to focus on the set of states in the tail of the distribution, the $\text{VaR}_{\alpha}$ is the smallest value in this set and, therefore, cannot detect the severity of the states with load shedding greater than $\text{VaR}_{\alpha}$, which is a drawback similar to the one associated to LOLP.

\subsubsection{$\text{CVaR}_{\alpha}$}
The \emph{Conditional Value-at-Risk} measures the expected value of the $\alpha$ worst states, i.e., the average of the states that comprise the tail of the probability distribution function of the load shedding and is defined as
\begin{equation}
	\label{eq:cvar_definition}
	\text{CVaR}_\alpha(R) = E[R : R \geq \text{VaR}_\alpha(R)]
\end{equation}

The $\text{CVaR}_{\alpha}$ index has several interesting properties \cite{AceTas2001expected} and \cite{RocUry2000optcvar} demonstrated that it is possible to formulate the $\text{CVaR}_{\alpha}$ calculation as a \emph{linear} programming problem independent of $\text{VaR}_{\alpha}$, which needs integer variables in its formulation.

As defined for the measures presented above, the SEP with $\text{CVaR}_\alpha$ constraint is shown in \eqref{eq:ProbPlan_InvOpeCnf_CVaR}.
\begin{subequations}
\label{eq:ProbPlan_InvOpeCnf_CVaR}
\begin{alignat}{5}
	\text{Min} \quad
	& \sum_{\igen \in \GenCSet} c_\igen x_\igen + O(x) & \\
	\text{s.t.:} \quad
	& b + \alpha^{-1} \sum_{s \in S} p^s y^s \leq \overline{\text{CVaR}} & \qquad \forall s \in S \label{eq:ProbPlan_InvOpeCnf_CVaR_maxcvar} \\
	& r^s \geq D - \sum_{\igen \in \GenSet} \xi_\igen^s \bar{g}_\igen x_\igen & \qquad \forall s \in S \\
	& r^s \geq 0 & \qquad \forall s \in S \\
	& y^s \geq  r^s - b & \qquad \forall s \in S \label{eq:ProbPlan_InvOpeCnf_CVaR_z} \\
	& y^s \geq 0 & \qquad \forall s \in S \\
	& x \in \mathcal{X}
\end{alignat}
\end{subequations}
where $\overline{\text{CVaR}}$ is a limit pre-established by the planner, $b$ is the variable that represents the implicitly calculated $\text{VaR}_{\alpha}$ while $y^s$ is the amount of load shedding that exceeds $b$, calculated by equation \eqref{eq:ProbPlan_InvOpeCnf_CVaR_z}. Therefore, the $\text{CVaR}_{\alpha}$ can be calculated as the sum of $b$ plus the expected value of $y^s$ conditioned to the states that exceed the $\text{VaR}_{\alpha}$, as shown in left-hand side of equation \eqref{eq:ProbPlan_InvOpeCnf_CVaR_maxcvar}.

\section{Considering Risk Constraints in Decomposition Schemes} \label{Decomposition}
As noted in previous sections, the reliability-constrained expansion planning problem is a large scale mixed-integer optimization problem. The incorporation of reliability constraints into the SEP requires the representation of a (possibly large) set of additional variables and constraints for \emph{every} state of the system and the number of states grows combinatorially with the number of generating units of the system.

The main objective of the design of mathematical decomposition techniques is to solve very complex or large problems through the repeated solution of a series of easier or smaller problems. From model \eqref{eq:ProbPlan_General} it is possible to note that the SEP with reliability constraints has a block structure and the problems are coupled by the investment decision vector $x$. This structure is opportune for the application of such techniques.

In this work, we use the Benders decomposition technique \cite{Ben1962PartProcMIP} to split the original problem into three sub-problems, intuitively reproducing the expansion planning process, which may consist of the following steps:
\begin{itemize}
\item First, the \emph{investment subproblem} (master) is solved, aiming for a new trial investment plan $x^{\mu}$, based on information obtained until iteration $\mu$: an approximation of the total cost function (investment plus approximate operation cost) and the approximation of the reliability ``feasible region'' (constraints representing the set of plans that meet the reliability criterion);
\item Given the proposed plan, $x^{\mu}$, the \emph{operation sub-problem} (slave) is solved and we check if the approximation of the cost function represented in the investment sub-problem is appropriate. If this function has not the adequate accuracy, a sensitivity analysis is conducted in order to build a new Benders optimality cut and improve the approximation of the cost function in the master problem;
\item For the same proposed plan $x^{\mu}$, the \emph{reliability sub-problem} (another slave) is solved to verify if the proposed solution is feasible with respect to the selected reliability index. If the solution is not feasible, a sensitivity analysis in this problem is carried out and a new Benders feasibility cut is obtained, improving the representation of the feasible region in the master problem.
\end{itemize}

In brief, at each iteration a trial solution is obtained from the master problem and sent to both slave sub-problems. Each sub-problem evaluates the decision $x^{\mu}$ in terms of its cost and reliability and also return Benders cuts for improving the representation of approximated operation cost and reliability index in the master problem. This procedure is repeated iteratively until a feasible solution with minimal total cost is found.

\subsection{Investment Subproblem}
The investment sub-problem can be formulated as the following mixed integer linear programming problem:
\begin{subequations}
\label{eq:ProbPlanDec_Inv}
\begin{alignat}{3}
	\text{Min} \quad
	& \sum_{\igen \in \GenSet} c_\igen x_{\igen} + \alpha & \\
	\text{s.t.:} \quad
	& \alpha \geq O(x^i) + \sum_{\igen \in \GenSet} \frac{\partial O(x^i)}{\partial x_\igen^i} (x_\igen - x_\igen^i) & \qquad i \in \mathcal{A} \label{eq:ProbPlanDec_Inv_OpeCuts} \\
	& R(x^i) + \sum_{\igen \in \GenSet} \frac{\partial R(x^i)}{\partial x_\igen^i} (x_\igen - x_\igen^i) \leq \bar{R} & \qquad i \in \mathcal{R} \label{eq:ProbPlanDec_Inv_RelCuts} \\
	& x \in \mathcal{X} &
\end{alignat}
\end{subequations}
where $\mathcal{A}$ and $\mathcal{R}$ are the sets of iterations where a cut has been added and $x^i$ is the trial solution vector found at iteration $i$. Constraints \eqref{eq:ProbPlanDec_Inv_OpeCuts} are called optimality cuts and are a first order approximation of the operation cost function $O(x)$. Similarly, constraints \eqref{eq:ProbPlanDec_Inv_RelCuts} are called feasibility cuts and are also a linear approximation of the feasibility region associated to the reliability index represented by the function $R(x)$.

\subsection{Operation sub-problem}
Given a trial investment decision $x^\mu$, the operation sub-problem can be formulated as
\begin{subequations}
\label{eq:ProbPlan_Ope_2}
\begin{alignat}{3}
	O(x^\mu) =
	\text{Min} \quad &
	\sum_{\igen \in \GenSet} d_\igen g_\igen + \defcost r & \\
	\text{s.t.:} \quad
	& \sum_{\igen \in \GenSet} g_\igen + r = D & \quad  \\
	& g_\igen \leq \bar{g}_{\igen} {x_\igen}^\mu & \quad \pi_\igen^{\bar{g}} & \quad \igen \in \GenSet \label{eq:ProbPlan_Ope_2_cappro}
\end{alignat}
\end{subequations}
where $\pi_\igen^{\bar{g}}$ corresponds to the dual variable of the maximum generation constraint of generator $\igen$.

From linear programming theory, it is known that $\pi_\igen^{\bar{g}}$ is the derivative of the objective function with respect to the right-hand side of constraints \eqref{eq:ProbPlan_Ope_2_cappro}. Applying the chain rule, we have
\begin{equation}
	\frac{\partial O(x^\mu)}{\partial {x_\igen}^\mu} = \frac{\partial O(x^\mu)}{\partial (\bar{g}_{\igen} {x_\igen}^\mu)} \times \frac{\partial (\bar{g}_{\igen} {x_\igen}^\mu)}{\partial {x_\igen}^\mu} = \pi_\igen^{\bar{g}} \times \bar{g}_{\igen}
\end{equation}
which is the derivative of the operation cost with respect to the investment decision $x^\mu$, used to construct the optimality Benders cuts \eqref{eq:ProbPlanDec_Inv_OpeCuts}.

In this section we considered the most simple single-stage deterministic representation of the system's operation. However, the construction of the Benders cut easily generalize to the multistage stochastic setting with all operative details. \cite{gorenstin1993power,campodonico2003expansion}

\subsection{Reliability sub-problem}
Given a trial investment decision $x^\mu$, it is possible to calculate the value of the risk measure associated to this plan, as well as the derivative of the reliability function $R(x^\mu)$ with respect to the investment decision, needed for the construction of the Benders cuts \eqref{eq:ProbPlanDec_Inv_RelCuts} that approximate the feasible region for the adopted reliability criterion in the master problem.

The feasibility Benders cuts consist in cutting planes that are tangent to the original feasibility region. Because of this reason, one of the requirements for the Benders decomposition method to be successfully applied is that the generated cuts cannot eliminate feasible solutions and this fact cannot be guaranteed if the sub-problem is non-convex. 

This same decomposition scheme cannot be implemented for SEP with LOLP and $\text{VaR}_{\alpha}$ constraints because they require the use of integer variables, characterizing a non-convex problem. It is still possible to use alternative and less efficient decompositions that remove infeasible integer solutions one by one. 
Therefore, this work considered only the SEP with EPNS or $\text{CVaR}_{\alpha}$ constraints as reliability criterion\footnote{Actually, it can be shown that is possible to approximate LOLP, EPNS and $\text{VaR}_{\alpha}$ criteria by just changing the parameters of $\text{CVaR}_{\alpha}$ planning model}.

\subsubsection{EPNS Criterion}
This sub-problem considers $R(x) = \text{EPNS}(x)$ and the solution can be obtained without the explicit representation of the load shedding variable $r_s$ because, given $x^\mu$, it is possible to calculate the EPNS as
\begin{equation}
\label{eq:ProbEENS_O}
	\text{EPNS}(x^\mu) = \sum_{s \in \Omega} p_s \times \left( D - \sum_{\igen \in \GenSet} \xi_{\igen s} \bar{g}_\igen {x_\igen}^\mu \right)
\end{equation}

Deriving the equation \eqref{eq:ProbEENS_O} with respect to the investment variable ${x_j}^\mu$, we have
\begin{equation}
\label{eq:ProbEENS_Derivada}
	\frac{\partial \text{EPNS}(x^\mu)}{\partial {x_j}^\mu} = -\sum_{s \in \Omega} p_s \xi_{js} \bar{g}_j \qquad \forall \igen \in \GenSet
\end{equation}
which consists in the coefficients of the feasibility Benders cuts \eqref{eq:ProbPlanDec_Inv_RelCuts}. For each investment plan, the $\text{EPNS}(x^\mu)$ and the derivative with respect to each ${x_\igen}^\mu$ can be obtained solving the states obtained from the Monte Carlo simulation independently. Then, a single cut is calculated by properly aggregating the results.

\subsubsection{$\text{CVaR}_{\alpha}$ Criterion}
The value for $\text{CVaR}_\alpha(x^\mu)$ if given by the solution of following linear programming sub-problem:
\begin{subequations}
\label{eq:ProbPlanDec_Cnf_CVaR}
\begin{alignat}{5}
    \text{Min} \quad
	& b + \alpha^{-1} \sum_{s \in S} p^s y^s \\
	\text{s.t.:} \quad
	& r^s \geq D - \sum_{\igen \in \GenSet} \xi_\igen^s \bar{g}_\igen {x_\igen}^\mu & \quad v^s & \quad \forall s \in S \label{eq:ProbPlanDec_Cnf_CVaR_r} \\
	& r^s \geq 0 & & \quad \forall s \in S \\
	& y^s \geq r^s - b & \quad w^s & \quad \forall s \in S \label{eq:ProbPlanDec_Cnf_CVaR_y} \\
	& y^s \geq 0 & & \quad \forall s \in S
\end{alignat}
\end{subequations}
where the decision variables $b$, $r^s$ and $y^s$ are non-negative and $v^s$ and $w^s$ are the dual variables of constraints \eqref{eq:ProbPlanDec_Cnf_CVaR_y} and \eqref{eq:ProbPlanDec_Cnf_CVaR_r}, respectively. The same solution process used for EPNS is applied for $\text{CVaR}_{\alpha}$ to obtain the feasibility cut.

\section{Case Study} \label{Study}
The case study consists of the expansion planning of the Bolivian (BO) generation system for a 7-years horizon. The system is composed by 28 existing generating hydro plants and 25 thermal plants with a total installed capacity of about 850MW (in the first year). In addition, 30 thermal projects are available as alternatives of investment for the expansion plan. Note that the maximum number of plants that may be operating in the system is 83 and, thus, the maximum number of states in the reliability sub-problem is $2^{83}$ (approximately  $10^{25}$ states).

Initially, we perform a comparative analysis for the expansion plans obtained with a hierarchical against an integrated methodology. Additionally, a comparison between the results obtained for the EPNS-constrained SEP and the $\text{CVaR}_{\alpha}$-constrained SEP is carried out.

We adopted a reliability limit for the EPNS equals to $1\%$ of the demand for each month over the horizon of study. The reliability sub-problem is solved by Monte Carlo simulation and it was considered that the convergence is reached when the coefficient of variation for the EPNS estimator is lower than $5\%$ \cite{PerMacOli1992AnaliticalAndMCPSAnal}.

The operation problem is solved by the SDDP algorithm in order to find the optimal dispatch under uncertainty, typically on inflow, demand, fuel costs, etc. All real-world details considered by the Bolivian ISO are represented. 

\subsection{Hierarchical Planning with EPNS criterion (HP-EPNS)}
\begin{figure}[!t]
	\centering
	\includegraphics[width=1.0\columnwidth]{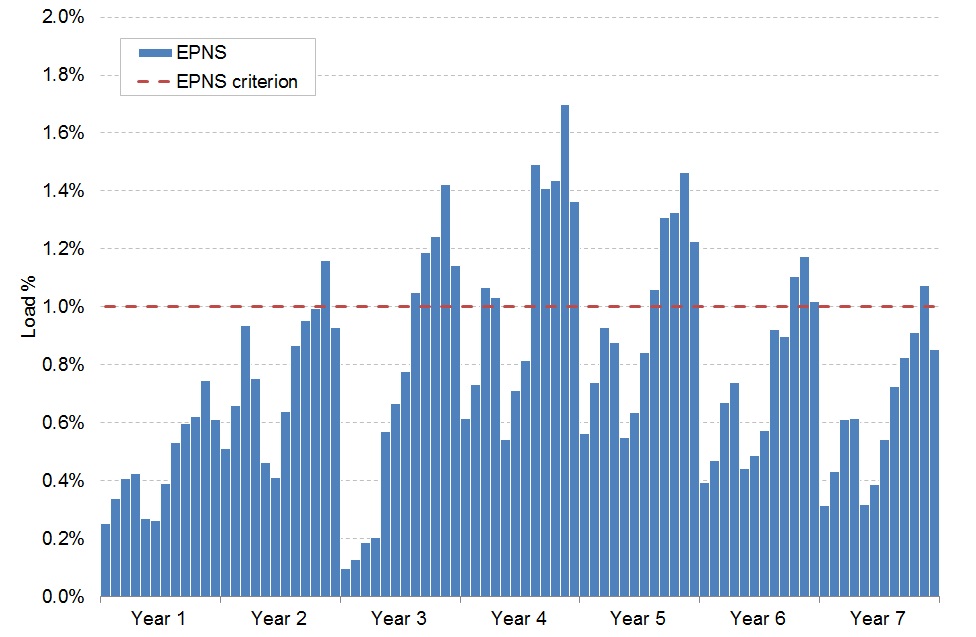}
 	\caption{Resulting EPNS for plan obtained with EP}
	\label{fig:epns_r02_epns}
\end{figure}

The hierarchical planning approach consists of obtaining an initial investment plan (first step) considering only the economic aspects, i.e., we aim to find the investment vector that minimizes the investment and operating costs (EP), as presented in section \ref{cap:PE}. Table \ref{tab:comp_plan_capadd} shows the additional capacity added in each year of the planning horizon and table \ref{tab:comp_plan_custo} its respective total cost. It can be seen in Figure \ref{fig:epns_r02_epns} that the resulting EPNS is not feasible according to the reliability criterion pre-established.

The second step of HP-EPNS is accomplished by considering the investment vector obtained in the first stage and solving the problem again to obtain reinforcements due reliability requirements. Observe in table \ref{tab:comp_plan_capadd} that it was necessary to invest in additional $44.1$MW in Year 2 (year that occurs the first violation of EPNS) to ensure feasibility for the reliability criterion. The associated total cost for this investment plan is also shown in table \ref{tab:comp_plan_custo}.

\begin{table}
\caption{Hierarchic x Integrated Planning: Additional Capacity [MW]}
\label{tab:comp_plan_capadd}
\begin{center}
\begin{tabular}{cccccccc}
\hline
Methodology & Y1 &            Y2 &            Y3 &   Y4 &    Y5 &   Y6 &            Y7 \\
\hline
         EP &  0 &    \textbf{0} &          65.8 & 44.1 & 285.6 & 44.1 &          44.1 \\
        HP-EPNS &  0 & \textbf{44.1} & \textbf{65.8} & 44.1 & 285.6 & 44.1 &          44.1 \\
        IP-EPNS &  0 &          44.1 & \textbf{44.1} & 44.1 & 285.6 & 44.1 & \textbf{44.1} \\
        IP-CVaR &  0 &          44.1 &          44.1 & 44.1 & 285.6 & 44.1 &    \textbf{0} \\
\hline
\end{tabular}
\end{center}
\end{table}

\subsection{Integrated Planning with EPNS criterion (IP-EPNS)}
In the integrated planning methodology both operation and reliability sub-problems are solved for each trial decision of the investment problem. Since the economic solution is not fixed when solving the reliability sub-problem as in the HP-EPNS, it is possible to consider in an integrated manner both the economic benefits and  the contribution for the attendance of the reliability criterion and then seek for the global optimal solution.

Table \ref{tab:comp_plan_capadd} shows that, in comparison with HP-EPNS, it is possible to invest in $44.1$MW in Year 3 instead of $65.8$MW and still get a reliability feasible plan. Note that compared to the methodology HP-EPNS, the IP-EPNS methodology obtains a lower cost investment plan but with higher operation cost. This fact illustrates that there is a benefit to invest in a generator with lower construction cost and higher operating cost because the reliability criterion could be met and total cost is lower, as shown in table \ref{tab:comp_plan_custo}.
\begin{table}[!t]
\caption{Hierarchical x Integrated Planning: Costs}
\label{tab:comp_plan_custo}
\begin{center}
\begin{tabular}{ccccc}
\hline
Methodology & Investment   & Operation   & Total       & \# of violated \\
            & Cost         & Cost        & Cost        & months         \\
\hline
EP          &    98.42 M\$ &  146.66 M\$ &  245.08 M\$ &    \textbf{22} \\
HP-EPNS     &   117.80 M\$ &  145.02 M\$ &  262.82 M\$ &              0 \\
IP-EPNS     &   100.06 M\$ &  152.17 M\$ &  252.23 M\$ &              0 \\
IP-CVaR     &    98.05 M\$ &  153.66 M\$ &  251.71 M\$ &              0 \\
\hline
\end{tabular}  
\end{center}
\end{table}

\subsection{Integrated Planning with $\text{CVaR}_{\alpha}$ criterion (IP-CVaR)}
In order to compare the risk measures EPNS and $\text{CVaR}_{\alpha}$, the $\text{CVaR}_{5\%}$ was calculated for the optimal investment plan obtained from the IP-EPNS and its maximum value (approximately 10\% of total load) was used as the $\text{CVaR}_{5\%}$ limit for the expansion planning.

\begin{figure}[!t]
	\centering
	\includegraphics[width=1.0\columnwidth]{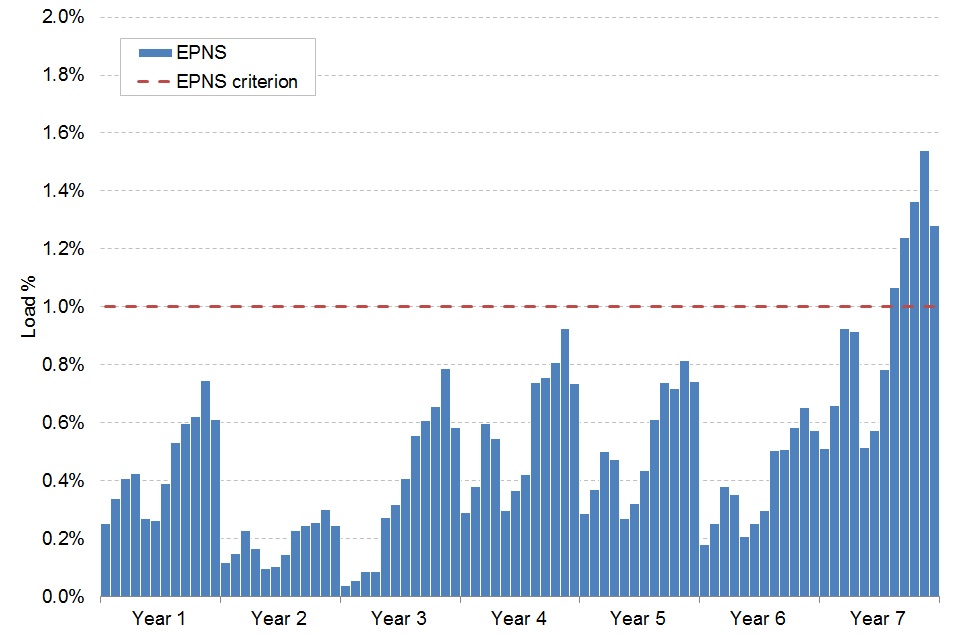}
 	\caption{Resulting EPNS for plan obtained with IP-CVaR}
	\label{fig:cvns_r06_epns}
\end{figure}

It is possible to observe the final costs obtained with this plan in table \ref{tab:comp_plan_custo} and the additional capacity in table \ref{tab:comp_plan_capadd}. The difference from the IP-EPNS is that the IP-CVaR methodology did not consider necessary to invest in $44.1$MW in the last year, obtaining an investment plan with lower cost. Furthermore, observing the EPNS for the plan in Figure \ref{fig:cvns_r06_epns}, note that, as the criterion of IP-CVaR methodology are related to the average of worst $5\%$ states and not the average of all states (EPNS), it was possible to find a more economic investment plan at a cost of EPNS exceeding $1\%$ of the load in the last year (but, as expected, the plan meets the $\text{CVaR}_{5\%} \leq 10\%$ of load criterion).

\section{Conclusions}\label{Conclusion}
This work presented a methodology to incorporate reliability constraints in the optimal power systems expansion planning. Through a real case example, it is shown that the simple application of an economic planning criterion is not enough to guarantee that the reliability criterion will be satisfied. It is also showed that the application of a hierarchical two-step procedure to solve the expansion planning problem with reliability criterion does not lead to the least cost solution.

It was shown that the integrated reliability constrained expansion planning problem can be modeled as a large scale MIP. In particular, we presented the MIP formulation for multiple reliability indexes. Benders decomposition can be used to exploit the problem structures and decouple the problem in a investment problem solved by standard MIP techniques; an operation problem solved SDDP; and a reliability problem solved by Monte Carlo simulation.

The advantage of using the integrated approach is to identify projects that contribute both in economic terms and in terms of improving overall system reliability, which might not possible by using the hierarchical procedure.

In addition to the traditional reliability measures commonly used in electrical systems, LOLP and EPNS, this work also illustrated how to incorporate the risk measures $\text{VaR}_{\alpha}$ and $\text{CVaR}_{\alpha}$, widely used in the financial area, into the power system expansion planning. It was showed that the use of the $\text{CVaR}_{\alpha}$ criterion allows the control of the depth of the probability distribution function of the system load shedding. This index can be a powerful tool of interest for system planners once it allows them to shape the expansion plan considering the prevention of disastrous events to the desired level of reliability by stimulating the diversification of the power plants ``portfolio''.

Moreover, this tool can be used to test and create policies for system expansion and to assess the individual contribution of each project in both economical and reliability aspects.

\section*{Acknowledgment}
The authors would like to acknowledge Nora Campod\'{o}nico and Silvio Binato from PSR and Marcia Fampa from COPPE/UFRJ for their valuable contributions and constant support.


\bibliographystyle{IEEEtran}
\bibliography{IEEEabrv,mylibrary}






\end{document}